\documentclass[a4paper,11pt]{article}
\usepackage{latexsym}
\usepackage{amssymb}
\usepackage{amsfonts}
\usepackage{amsmath}
\usepackage{indentfirst}
\usepackage[dvips]{graphicx}
\usepackage{epsfig}
\usepackage{rotating}
\newtheorem{prop}{Proposition}[section]
\newtheorem{teor}{Theorem}[section]
\newtheorem{lemma}{Lemma}[section]
\newtheorem{cor}{Corollary}[section]
\newcommand{\nkinN}{n,k\in \mathbf{N}}
\newcommand{\ninN}{n\in \mathbf{N}}

\newcommand{\kinN}{k\in \mathbf{N}}
\newcommand{\cvd}{\quad $\blacksquare$\bigskip}

\date{}

\author{Silvia Bacchelli\thanks{I.I.S. ``Crescenzi-Pacinotti", via Saragozza 9,
40123 Bologna, Italy \tt{silvia.bacchelli@istruzione.it}}\and Luca
Ferrari\thanks{Dipartimento di Sistemi e Informatica, viale
Morgagni 65, 50134 Firenze, Italy \tt{ferrari@dsi.unifi.it,
pinzani@dsi.unifi.it, renzo.sprugnoli@unifi.it}}\and Renzo
Pinzani$^\dag$\and Renzo Sprugnoli$^\dag$}

\title{Mixed succession rules: the commutative case} 

\begin{document}

\maketitle

\begin{abstract}
We begin a systematic study of the enumerative combinatorics of
\emph{mixed succession rules}, which are succession rules such
that, in the associated generating tree, the nodes are allowed to
produce their sons at several different levels according to
different production rules. Here we deal with a specific case,
namely that of two different production rules whose rule operators
commute. In this situation, we are able to give a general formula
expressing the sequence associated with the mixed succession rules
in terms of the sequences associated with the component production
rules. We end by providing some examples illustrating our
approach.
\end{abstract}

\section{Introduction}\label{introduzione}

Among the many methods that have been developed to enumerate
combinatorial structures, the role of the \emph{ECO method} has
been growing in the last decade, thanks to its intrinsic
simplicity and to the effectiveness of the combinatorial
constructions it generates. The variety of problems in which the
ECO method has shown its soundness ranges from enumerative and
bijective combinatorics to random and exhaustive generation.

The roots of the ECO method can be traced back to \cite{CGHK},
where the authors study Baxter permutations and introduce for the
first time the concept of a \emph{generating tree}. Successively,
West \cite{W1,W2} introduced the notion of a \emph{succession
rule} to give a formal description of generating trees in the
context of permutation enumeration and Barcucci, Del Lungo,
Pergola and Pinzani \cite{BDLPP} extended the technique of
generating trees, finding a general way of constructing
combinatorial objects which can be often described using such
formal tools.

\bigskip

The \emph{classical ECO method} (a detailed description of which
can be found, for instance, in \cite{BDLPP}) consists of a
recursive construction for a class of objects by means of an
operator which performs a ``local expansion" on the objects
themselves. Typically, starting from an object \emph{of size $n$},
an ECO construction allows to produce a set of new objects,
\emph{of size $n+1$}, in such a way that, iterating the
construction, all the objects of the class are obtained precisely
once. If the construction is sufficiently regular, it can be often
described by means of a \emph{succession rule}, which is a system
of the form
\begin{equation}\label{Omega}
\Omega :\left\{ \begin{array}{ll} (a)
\\ (k)\rightsquigarrow (e_1 (k))(e_2 (k))\cdots (e_k (k))
\end{array}\right. .
\end{equation}

The meaning of what is written above is the following. Each object
of the class is given a label $(k)$. When performing the ECO
construction, an object labelled $(k)$ produces $k$ new objects
labelled, respectively, $(e_1 (k)),(e_2 (k)),\ldots ,(e_k (k))$.
Moreover, the object of minimum size has label $(a)$ (which is
called the \emph{axiom} of the succession rule). To have a
graphical description of a succession rule one usually draws its
\emph{generating tree}, that is the infinite, rooted, labelled
tree whose root is labelled $(a)$ (like the axiom) and such that
each node labelled $(k)$ has $k$ sons, labelled $(e_1 (k)),(e_2
(k)),\ldots ,(e_k (k))$, respectively. It is evident from this
definition that we can introduce a notion of \emph{level} on
generating trees, by saying that the root lies at level 0, and a
node lies at level $n$ when its parent lies at level $n-1$.

We remark that, from the above given definition, a node labelled
$(k)$ has precisely $k$ sons. When a succession rule has this
property it is often said to be \emph{consistent}. However, one
can also consider succession rules (and generating trees) in which
the value of a label does not necessarily represent the number of
its sons, and this will be frequently done in the sequel.
Moreover, we would like to warn the reader that, even if we will
sometimes give definitions using consistent succession rules
(since this is the convention when working with the ECO method),
we will constantly make use of succession rules which are not
necessarily consistent.

From the enumerative point of view, the main information encoded
in a generating tree (and thus in its associated succession rule)
is given by the \emph{level polynomial} $p_n
(x)=\sum_{k}p_{n,k}x^k$, defined by setting $p_{n,k}$ equal to the
number of nodes labelled $k$ at level $n$, and by the
\emph{associated integer sequence} $(f_n )_{\ninN}$, which is
defined, in terms of the level polynomials, as $f_n =p_n (1)$, and
represents the total number of nodes at level $n$.

We point out that the infinite lower triangular array
$(p_{n,k})_{\nkinN}$, sometimes called the \emph{AGT matrix}
\cite{MV}, or \emph{ECO matrix} \cite{FP2}, often happens to be a
\emph{Riordan array} \cite{Sp}. By definition, this means that
every element $p_{n,k}$ can be expressed by using a pair of formal
power series $(d(t),h(t))$ in such a way that it is precisely the
coefficient of $t^n$ of $d(t)h(t)^k$. In this case, many counting
properties of the generating tree can be found in an algebraic
way, by using the related theory.

To give a quick example, consider the succession rule
\begin{displaymath}
\Omega:  \left\{ \begin{array}{ll}
(2)\\
(k)\rightsquigarrow (2)(3)(4)\cdots (k)(k+1)
\end{array} \right. ,
\end{displaymath}
defining Catalan numbers $1,2,5,14,42,132,\ldots$ (see, for
example, \cite{BDLPP}). The first levels of its generating tree
can be depicted as follows:

\begin{center}
\setlength{\unitlength}{0.8mm}
\begin{picture}(100,100)
\put(50,90){\makebox(0,0){2}} \put(25,70){\makebox(0,0){2}}
\put(75,70){\makebox(0,0){3}} \put(10,50){\makebox(0,0){2}}
\put(40,50){\makebox(0,0){3}} \put(60,50){\makebox(0,0){2}}
\put(75,50){\makebox(0,0){3}} \put(90,50){\makebox(0,0){4}}
\put(5,30){\makebox(0,0){2}} \put(15,30){\makebox(0,0){3}}
\put(35,30){\makebox(0,0){2}} \put(40,30){\makebox(0,0){3}}
\put(45,30){\makebox(0,0){4}} \put(55,30){\makebox(0,0){2}}
\put(65,30){\makebox(0,0){3}} \put(70,30){\makebox(0,0){2}}
\put(75,30){\makebox(0,0){3}} \put(80,30){\makebox(0,0){4}}
\put(85,30){\makebox(0,0){2}} \put(90,30){\makebox(0,0){3}}
\put(95,30){\makebox(0,0){4}} \put(100,30){\makebox(0,0){5}}
\put(27,72){\line(5,4){20}} \put(73,72){\line(-5,4){20}}
\put(12,52){\line(3,4){11}} \put(38,52){\line(-3,4){11}}
\put(62,52){\line(3,4){11}} \put(75,53){\line(0,1){14}}
\put(88,52){\line(-3,4){11}} \put(6,32){\line(1,4){3.5}}
\put(14,32){\line(-1,4){3.5}} \put(36,32){\line(1,4){3.5}}
\put(40,32){\line(0,1){14}} \put(44,32){\line(-1,4){3.5}}
\put(56,32){\line(1,4){3.5}} \put(64,32){\line(-1,4){3.5}}
\put(71,32){\line(1,4){3.5}} \put(75,32){\line(0,1){14}}
\put(79,32){\line(-1,4){3.5}} \put(86,32){\line(1,4){3.5}}
\put(90,32){\line(0,1){14}} \put(94,32){\line(-1,4){3.5}}
\put(98,32){\line(-1,2){7}}
\end{picture}
\end{center}

\vspace{-15mm}

Here the level polynomials are $p_n (x)=\sum_{k}b_{n,k}x^k$, where
$b_{n,k}=\frac{k-1}{2n-k+3}{2n-k+3\choose n+1}$ are the usual
ballot numbers. Notice that, in this example, by shifting the
column index $k$ by two positions, it can be shown that
$B=(b_{n,k})_{\nkinN}=(C(t),tC(t))$, where
$C(x)=\frac{1-\sqrt{1-4x}}{2x}$ is the generating function of
Catalan numbers. Using the results of \cite{Sp}, we can compute
the row sums and the weighted row sums of $B$, thus obtaining the
basic data for evaluating the probability distribution of the
labels in the generating tree.

\bigskip

The ECO method has been fruitfully applied to several problems,
not only of an enumerative nature. For instance, using this
technique it has been possible to develop efficient algorithms for
the random \cite{BDLP} and exhaustive \cite{BBGP} generation of
combinatorial objects. Moreover, in \cite{BGPP} the authors
describe a general exhaustive generation algorithm (working for a
wide family of structures) defining Gray codes not depending on
the specific nature of the objects to be generated, but only on
the properties of a succession rule encoding an ECO construction
for the objects under consideration. Due to its plentiful
applications, it is then worth exploring in a deeper way the
features and capabilities of such a method.

\bigskip

Despite its wide range of applicability, there are many
combinatorial constructions which cannot be naturally described by
using the classical ECO method (and classical succession rules) as
exploited above. For instance, in \cite{BM} a generalization of
the method is considered, allowing succession rules in which the
labels are pairs of integers (rather than integers). Applications
of this generalized method to the enumeration of pattern avoiding
permutations are shown in the cited paper.

A very strong limitation in the possibility of describing a
combinatorial construction by means of a succession rule lies in
the definition we have given of this fundamental tool. The
generating tree of a succession rule has the property that, if the
level of a node is $n$, then the level of all its sons is $n+1$.
From a combinatorial point of view, this means that a (classical)
ECO construction performed on an object of a given size produces
objects of the successive size. However, it may well happen that a
combinatorial construction, having all the reasonable features to
be called ECO, does not behave in the standard way with respect to
the notion of size. More precisely, starting from an object of
size $n$, we can construct new objects whose sizes are greater
than $n$ (but not necessarily equal to $n+1$). The formalization
of these concepts leads to the notion of what can be called a
\emph{mixed succession rule}. Roughly speaking, the idea is to
consider a set of (possibly different) succession rules acting on
the objects of a class and producing sons at different levels. To
be more formal, we introduce here the simplest instance of this
general situation, by considering two succession rules producing
their sons at the two successive levels. These will be called
\emph{doubled mixed succession rules}. Given two succession rules
$\Omega$ as in (\ref{Omega}) and
\begin{displaymath}
\Sigma :\left\{ \begin{array}{ll} (b)
\\ (k)\rightsquigarrow (d_1 (k))(d_2 (k))\cdots (d_k (k))
\end{array}\right. ,
\end{displaymath}
we define the \emph{doubled mixed succession rule associated with
the pair $(\Omega ,\Sigma )$ with axiom $(c)$} to be the
succession rule $(c)\Omega^{+1}\Sigma^{+2}$, defined by
\begin{displaymath}
(c)\Omega^{+1}\Sigma^{+2} :\left\{ \begin{array}{ll} (c)
\\ (k)&\stackrel{+1}{\rightsquigarrow}(e_1 (k))(e_2 (k))\cdots (e_k (k))
\\ &\stackrel{+2}{\rightsquigarrow}(d_1 (k))(d_2 (k))\cdots (d_k (k))
\end{array}\right. .
\end{displaymath}

The generating tree associated with $(c)\Omega^{+1}\Sigma^{+2}$
has the property that each node labelled $(k)$ lying at level $n$
produces two sets of sons, the first set being $(e_1 (k)),(e_2
(k)),\ldots ,(e_k (k))$ at level $n+1$ and the second one being
$(d_1 (k)),(d_2 (k)),\ldots ,(d_k (k))$ at level $n+2$ (so that it
produces a total of $2k$ sons).

To justify our interest in this kind of notion, we remark that
instances of (general) mixed succession rules have occasionally
surfaced in some previous works; to cite only one example, in
\cite{GPP} vexillary involutions are enumerated by making use of a
specific mixed succession rule. The first systematic treatment of
mixed succession rules has been undertaken in \cite{FPPR2}, where
the special case $\Sigma =\Omega$ has been examined in great
detail. The present paper represents the first attempt to tackle
the general case, aiming at developing a general theory of mixed
succession rules. More precisely, the main goal would be to
succeed in expressing the sequence associated with a mixed
succession rule $(c)\Omega^{+1}\Sigma^{+2}$ in terms of the
sequences associated with $\Omega$ and $\Sigma$, possibly changing
the axioms. For this reason, in section \ref{prodrule} we study
some enumerative properties of what we have called
\emph{production rules}, which are, by definition, succession
rules without the axiom. The problem of studying mixed succession
rules in its full generality seems quite difficult; it is somehow
related to the theory of power series in several
\emph{noncommuting} variables. Here we deal only with a special
case, namely when the two rule operators (see section
\ref{preliminari}) of $\Omega$ and $\Sigma$ commute. In this
situation, we are able to find a general formula for the sequence
associated with $(c)\Omega^{+1}\Sigma^{+2}$; moreover, we also
describe some examples our theory can be applied to.

\bigskip

We would like to remark that this problem has also been considered
from the point of view of Riordan arrays \cite{BMS}. Each Riordan
array determines a specific sequence $(a_k)_{\kinN}$, called the
\emph{$A$-sequence} of the array, such that, for every $\nkinN$:
\begin{displaymath}
b_{n+1,k+1}=\sum_{k=0}^{\infty}a_k b_{n,n+k}.
\end{displaymath}

When the $A$-sequence contains integer numbers only, it is related
to the succession rule of an associated generating tree (if any)
as shown in \cite{MV}. However, it can happen that the
$A$-sequence has a complicated expression, whereas the $A$-matrix
(as defined in \cite{BMS}) is simple. This corresponds to an ECO
construction in which the elements of size $n$ produce objects of
different sizes (greater than $n+1$).

\bigskip

In closing this introduction, we recall some notations we will
frequently use in the next pages.

The sets of natural and real numbers will be denoted $\mathbf{N}$
and $\mathbf{R}$, respectively.

The following linear operators on the vector space of one-variable
polynomials will be often considered: $\mathbf{x}$ (respectively,
$\mathbf{t}$) is the operator of multiplication by $x$
(respectively, $t$), $D$ is the usual derivative operator, and $T$
is the \emph{factorial derivative operator}, which is, by
definition, the linear operator mapping $x^n$ into $1+x+\cdots
+x^{n-1}=\sum_{i=0}^{n-1}x^i$ (for $n\geq 1$) and 1 into 0.

\section{Preliminaries on rule operators}\label{preliminari}

Given a succession rule $\Omega$ as in (\ref{Omega}), we can
associate with it a linear operator on the vector space of
one-variable polynomials $K[x]$, to be denoted $L=L_{\Omega}$ (the
subscript will be omitted when it is clear from the context). To
define such an operator, we use the canonical basis $(x^n
)_{\ninN}$:
\begin{eqnarray*}
L&:&K[x]\longrightarrow K[x]
\\&:&1\longmapsto x^a
\\&:&x^k \longmapsto x^{e_1 (k)}+\cdots +x^{e_k (k)},\qquad
\textnormal{if $k$ appears in $\Omega$,}
\\&:&x^h \longmapsto hx^h ,\qquad \textnormal{otherwise.}
\end{eqnarray*}

As it is easy to understand, each enumerative property of $\Omega$
can be suitably translated into some property of $L$. For
instance, if $(f_n )_{\ninN}$ is the numerical sequence associated
with $\Omega$, for any $\ninN$ the following equality holds:
\begin{displaymath}
f_n =[L^{n+1}(1)]_{x=1},
\end{displaymath}
where we have used square brackets to denote the operator of
evaluation at a specific value. The linear operator $L$ is called
the \emph{rule operator} associated with $\Omega$. We refer the
reader to \cite{FPPR1,FPPR3,FP1,FP2} for the definition, the main
properties and some applications of this notion.

\bigskip

\emph{Remarks.}
\begin{enumerate}
\item In the above definition we have denoted with $K$ a generic
fields of coefficients. For what concerns us here, it is largely
enough to take $K=\mathbf{R}$. Actually, all the theory of rule
operators could be equally developed on the semiring module
$\mathbf{N}[x]$ of polynomials with nonnegative integer
coefficients. \item If $p_n (x)=\sum_{k}p_{n,k}x^k$ is the $n$-th
level polynomial of $\Omega$, it is
\begin{displaymath}
L_{\Omega}(p_n (x))=p_{n+1}(x).
\end{displaymath}
\end{enumerate}

\bigskip

\emph{Examples.} Here are some examples of rule operators
associated with more or less well-known succession rules. All
these examples can be found, for example, in \cite{FP1}.
\begin{enumerate}
\item A succession rule for factorial numbers:
\begin{eqnarray*}
\left\{ \begin{array}{ll} (1)
\\ (k)\rightsquigarrow (k+1)^k
\end{array}\right. ;
\end{eqnarray*}
associated rule operator:
\begin{eqnarray*}
L(1)&=&x
\\ L(x^k )&=&(\mathbf{x}^2 D)(x^k )=kx^{k+1},\qquad k\geq 1.
\end{eqnarray*}
\item A succession rule for arrangements:
\begin{eqnarray*}
\left\{ \begin{array}{ll} (2)
\\ (k)\rightsquigarrow (k)(k+1)^{k-1}
\end{array}\right. ;
\end{eqnarray*}
associated rule operator:
\begin{eqnarray*}
L(1)&=&x^2
\\ L(x^k )&=&(\mathbf{x}^2 D-\mathbf{x}+1)(x^k )=x^k +(k-1)x^{k+1},\qquad k\geq 1.
\end{eqnarray*}
\item A succession rule for involutions:
\begin{eqnarray*}
\left\{ \begin{array}{ll} (1)
\\ (k)\rightsquigarrow (k-1)^{k-1}(k+1)
\end{array}\right. ;
\end{eqnarray*}
associated rule operator:
\begin{eqnarray*}
L(1)&=&x
\\ L(x^k )&=&(D+\mathbf{x}-\mathbf{x}^{-1})(x^k )=(k-1)x^{k-1}+x^{k+1},\qquad k\geq 1.
\end{eqnarray*}
\item A succession rule for Bell numbers:
\begin{eqnarray*}
\left\{ \begin{array}{ll} (1)
\\ (k)\rightsquigarrow (k)^{k-1}(k+1)
\end{array}\right. ;
\end{eqnarray*}
associated rule operator:
\begin{eqnarray}\label{bell}
L(1)&=&x\nonumber
\\ L(x^k )&=&(\mathbf{x}D+\mathbf{x}-1)(x^k )=(k-1)x^k +x^{k+1},\qquad k\geq 1.
\end{eqnarray}
\end{enumerate}

Further examples involving the factorial derivative operator $T$
can be considered.

\begin{itemize}
\item[5.] A succession rule for Catalan numbers:
\begin{eqnarray*}
\left\{ \begin{array}{ll} (1)
\\ (k)\rightsquigarrow (2)(3)\cdots (k)(k+1)
\end{array}\right. ;
\end{eqnarray*}
associated rule operator:
\begin{eqnarray}\label{catalan}
L(1)&=&x\nonumber
\\ L(x^k )&=&(\mathbf{x}^2 T)(x^k )=\sum_{i=2}^{k+1}x^i,\qquad k\geq 1.
\end{eqnarray}
\item[6.] A succession rule for Motzkin numbers:
\begin{eqnarray*}
\left\{
\begin{array}{ll} (1)
\\ (1)\rightsquigarrow (2)
\\ (k)\rightsquigarrow (1)(2)\cdots (k-1)(k+1)\qquad k\geq 1
\end{array}\right. ;
\end{eqnarray*}
associated rule operator:
\begin{eqnarray}\label{motzkin}
L(1)&=&x\nonumber
\\ L(x^k )&=&(\mathbf{x}T+\mathbf{x}-\mathbf{1})(x^k )=\sum_{i=1}^{k-1}x^i +x^{k+1},\qquad k\geq 2.
\end{eqnarray}
\end{itemize}



\section{Production rules}\label{prodrule}

With the expression \emph{production rule} we will mean here a
succession rule without its axiom. Hence the generic form of a
production rule is
\begin{equation}\label{produzione}
(k)\rightsquigarrow (e_1 (k))(e_2 (k))\cdots (e_k (k)).
\end{equation}

Clearly, in the same way a succession rule determines a unique
numerical sequence, a production rule defines a family of
sequences $(f_n ^{(a)})_{\ninN}$, depending on the axiom $(a)$
which we choose for the rule (\ref{produzione}).

\bigskip

From now on, given a succession rule $\Omega$ as in (\ref{Omega}),
we will denote with $L_a$ the associated rule operator, $(a)$
being the axiom of $\Omega$. Using this terminology, given a
production rule as in (\ref{produzione}), the family of operators
$(L_a )_{a\in \mathbf{N}}$ will be called the \emph{family of rule
operators} associated with the production rule. In this section we
will be interested in finding formulas to relate the various
sequences associated with the same production rule. To this aim,
we start by observing the following (very easy but quite
important) facts.

\begin{enumerate}
\item For any $a,b\in \mathbf{N}$, we have $L_a (x^k )=L_b (x^k
)$, when $k\neq 0$, and $L_a (1)=x^a$, $L_b (1)=x^b$. For this
reason, in what follows, for $k\neq 0$, we will simply write
$L(x^k )$, without specifying the axiom, and we will speak of
\emph{the} rule operator associated with the production rule
whenever we restrict our attention to the subspace $xK[x]$ (i.e.
the subspace spanned by the positive powers of $x$). \item The
$n$-th term of the numerical sequence of the family with axiom
$(b)$, that is $f_n ^{(b)}$, can be computed using the following
formula:
\begin{displaymath}
f_n ^{(b)}=[L^n (x^b )]_{x=1}.
\end{displaymath}
In the same way, to compute $f_n ^{(b+1)}$ we get:
\begin{displaymath}
f_n ^{(b+1)}=[L^n (x^{b+1})]_{x=1}=[L^n \mathbf{x}(x^b )]_{x=1}.
\end{displaymath}
Then, if one knows the \emph{Pincherle derivative} of $L$, which
is, by definition, the operator $L'=L\mathbf{x}-\mathbf{x}L$ (see
\cite{RKO}), it should be possible to express the operator $L^n
\mathbf{x}$ as a linear combination of monomials of the kind
$\mathbf{x}^{\alpha}L^{\beta}$. This should allow, at least in
principle, to obtain an expression for$f_n ^{(b+1)}$ in terms of
known quantities (namely $f_m ^{(b+1)}$, with $m<n$, and $f_k
^{(a)}$, with $a\leq b$). To better understand how to proceed in
the concrete cases, it is convenient to have a look at a specific
example.
\end{enumerate}

\subsection{A Bell-like production rule}

Let us consider the following production rule:
\begin{displaymath}
\omega :\quad (k)\rightsquigarrow (k)^{k-1}(k+1).
\end{displaymath}

On the nonzero powers of $x$, the rule operator associated with
$\omega$ acts as follows:
\begin{displaymath}
L(x^k )=(\mathbf{x}D+\mathbf{x}-\mathbf{1})(x^k ),\qquad \forall
k\neq 0.
\end{displaymath}

Thus, on $xK[x]$ we have that
$L=\mathbf{x}D+\mathbf{x}-\mathbf{1}$. Concerning the Pincherle
derivative $L'$, on $xK[x]$ we have:
\begin{eqnarray*}
&&L\mathbf{x}=(\mathbf{x}D+\mathbf{x}-\mathbf{1})\mathbf{x}=\mathbf{x}(\mathbf{x}D+\mathbf{1})+\mathbf{x}^2
-\mathbf{x}=\mathbf{x}^2 D+\mathbf{x}^2 ,
\\ &&\mathbf{x}L=\mathbf{x}(\mathbf{x}D+\mathbf{x}-\mathbf{1})=\mathbf{x}^2
D+\mathbf{x}^2 -\mathbf{x},
\end{eqnarray*}
(here we have used the well-known identity
$D'=D\mathbf{x}-\mathbf{x}D=\mathbf{1}$). Therefore, we get:
\begin{equation}\label{pincherle}
L'=\mathbf{x}.
\end{equation}

Our next goal is to express the operator $L^n \mathbf{x}$ in such
a way that some identity between the terms of the sequences
associated with $\omega$ can be determined.

\begin{prop} For any $n\geq 1$, we have:
\begin{equation}\label{esempio}
L^n \mathbf{x}=\sum_{k=0}^{n-1}{n-1 \choose
k}\mathbf{x}L^{k+1}+L^{n-1}\mathbf{x}.
\end{equation}
\end{prop}

\emph{Proof.} First observe that, from (\ref{pincherle}), we get
$L\mathbf{x}=\mathbf{x}(\mathbf{1}+L)$, whence:
\begin{displaymath}
\mathbf{x}(\mathbf{1}+L)^n
=\mathbf{x}(\mathbf{1}+L)(\mathbf{1}+L)^{n-1}=L\mathbf{x}(\mathbf{1}+L)^{n-1}.
\end{displaymath}
and then, by iterating:
\begin{displaymath}
\mathbf{x}(\mathbf{1}+L)^n =L^n \mathbf{x}
\end{displaymath}

The above formula is easily seen to hold for $n\geq 0$. By
iterating, we then get:
\begin{eqnarray*}
L^n \mathbf{x}&=&\mathbf{x}(\mathbf{1}+L)^n
=\mathbf{x}(\mathbf{1}+L)^{n-1}+\mathbf{x}L(\mathbf{1}+L)^{n-1}
\\ &=&L^{n-1}\mathbf{x}+\mathbf{x}L\sum_{k=0}^{n-1}{n-1 \choose
k}L^{k}
\\ &=&L^{n-1}\mathbf{x}+\sum_{k=0}^{n-1}{n-1 \choose
k}\mathbf{x}L^{k+1}.\quad \blacksquare
\end{eqnarray*}

\bigskip


Formula (\ref{esempio}) yields an interesting recursive expression
for the family of sequences associated with $L$.

\begin{cor} Denoting by $(f_n ^{(b)})_{\ninN}$ the sequence
associated with $L$ with axiom $(b)$, it is
\begin{displaymath}
f_n ^{(b+1)}=\sum_{k=0}^{n-1}{n-1\choose
k}f_{k+1}^{(b)}+f_{n-1}^{(b+1)}.
\end{displaymath}
\end{cor}

\emph{Proof.}\quad Indeed, using the above proposition, we get:
\begin{eqnarray*}
f_n ^{(b+1)}&=&[L^n \mathbf{x}(x^b
)]_{x=1}=\sum_{k=0}^{n-1}{n-1\choose k}[\mathbf{x}L^{k+1}(x^b
)]_{x=1}+[L^{n-1}\mathbf{x}(x^b )]_{x=1}
\\ &=&\sum_{k=0}^{n-1}{n-1\choose k}f_{k+1}^{(b)}+f_{n-1}^{(b+1)},
\end{eqnarray*}
as desired.\cvd

As a further consequence, one can, for example, determine the
generating function $f^{(b)}(x)$ of $(f_n ^{(b)})_{\ninN}$.
Standard generating function arguments provide the following
result.

\begin{cor} Denoting by $\mathcal{B}$ the binomial transform
operator on exponential generating functions (such as, given an
exponential generating function $f(x)$, it is
$\mathcal{B}(f(x))=e^x f(x)$), we have:
\begin{displaymath}
f^{(b)}(x)=\mathcal{B}^{b-1}(f^{(1)}(x)).
\end{displaymath}
Since $f^{(1)}(x)=e^{e^{x}-1}$ is the exponential generating
function of Bell numbers $1,1,2,5,15,52,203,\ldots$, we get
\begin{displaymath}
f^{(b)}(x)=e^{e^{x}+(b-1)x-1}
\end{displaymath}
\end{cor}

\subsection{Further examples}\label{further}

It is obvious that what we have done for the rule operator
$\mathbf{x}D+\mathbf{x}-\mathbf{1}$ can be done (at least in
principle) for any other rule operator. In particular, we mention
here two further cases, giving only the final results, and leaving
to the reader the details of the proofs.

\begin{description}
\item[1)] \emph{A Catalan-like production rule.}\quad Consider the
production rule
\begin{displaymath}
\omega :\quad (k)\rightsquigarrow (2)(3)\cdots (k)(k+1),
\end{displaymath}
which is related to Catalan numbers. Indeed, it is known that the
generating function of the sequence with axiom $(b)$ is $C^b (x)$
(recall that $C(x)$ is the generating function of Catalan
numbers). We can use our approach to rediscover this result.

Indeed, we first observe that $L=\mathbf{x}^2 T$ is the rule
operator associated with $\omega$ (as usual, we restrict our
attention to the subspace $xK[x]$). We start by computing the
Pincherle derivative $L'$.

\begin{lemma} Let $c_2 :xK[x]\longrightarrow xK[x]$ be the linear
operator defined on the canonical basis by setting $c_2 (x^n
)=x^2$, for every $n>0$. Then we have that
\begin{displaymath}
L'=c_2 .
\end{displaymath}
\end{lemma}

From this lemma we obtain the following result, which is essential
in deriving the final recursion.

\begin{prop} For every $\ninN$, it holds:
\begin{displaymath}
L^n \mathbf{x}=\mathbf{x}L^n +\sum_{i=0}^{n-1}L^i c_2 L^{n-1-i}.
\end{displaymath}
\end{prop}

Finally, after having observed that, for all $i\in \mathbf{N}$, it
is clearly $[L^i c_2 (p(x))]_{x=1}=f_i ^{(2)}\cdot p(1)$, we get
the following recursion for the family of sequences associated
with the starting production rule.

\begin{prop} For every $b,\ninN$, we have:
\begin{displaymath}
f_n ^{(b+1)}=f_n ^{(b)}+\sum_{i=0}^{n-1}f_i ^{(2)}f_{n-1-i}^{(b)}.
\end{displaymath}
As a consequence, the generating function of $(f_n
^{(b)})_{\ninN}$ is given by the convolution $C^b (x)$ of the
generating function $C(x)=\frac{1-\sqrt{1-4x}}{2x}$ of Catalan
numbers.
\end{prop}

\item[2)] \emph{A Motzkin-like production rule.}\quad Consider the
production rule
\begin{displaymath}
\omega :\quad (k)\rightsquigarrow (1)(2)\cdots (k-1)(k+1),
\end{displaymath}
which is related to Motzkin numbers. Namely, it is known that
$\omega$, with axiom $(1)$, yields the sequence of Motzkin numbers
$1,1,2,4,9,21,\ldots$. The rule operator associated with $\omega$
is $L=\mathbf{x}T+\mathbf{x}-1$ (also in this case we refer to
$xK[x]$). We can proceed analogously as we have done in the
previous example.

\begin{lemma}Let $c_1 :xK[x]\longrightarrow xK[x]$ be the linear
operator defined on the canonical basis by setting $c_1 (x^n )=x$,
for every $n>0$. Then we have that
\begin{displaymath}
L'=c_1 .
\end{displaymath}
\end{lemma}

\begin{prop} For every $\ninN$, it holds:
\begin{displaymath}
L^n \mathbf{x}=\mathbf{x}L^n +\sum_{i=0}^{n-1}L^i c_1 L^{n-1-i}.
\end{displaymath}
\end{prop}

\begin{prop} For every $b,\ninN$, we have:
\begin{displaymath}
f_n ^{(b+1)}=f_n ^{(b)}+\sum_{i=0}^{n-1}f_i ^{(1)}f_{n-1-i}^{(b)}.
\end{displaymath}
As a consequence, the generating function of $(f_n
^{(b)})_{\ninN}$ is given by the convolution
$f^{(b)}(x)=f^{(b-1)}(x)(xM(x)+1)$, where
$M(x)=\frac{1-x-\sqrt{1-2x-3x^2 }}{2x^2}$ is the generating
function of Motzkin numbers, and so $f^{(b)}(x)=M(x)\cdot
(xM(x)+1)^{b-1}$.
\end{prop}

\end{description}

\section{Commuting rule operators}

In this section we finally tackle the problem of finding an
explicit expression for the numerical sequence associated with a
mixed succession rule. As we have said in the introduction, in
what follows we will deal with the simplest case of a mixed
succession rule, namely the case of a \emph{doubled rule}: this
means that, in the associated generating tree, each node at level
$n$ produces a set of sons at level $n+1$, according to a
production rule $\omega$, and another set of sons at level $n+2$,
according to another production rule $\sigma$. If $(b)$ is the
axiom of the doubled rule, such a generating tree can be
synthetically represented as follows:

\begin{figure}[!h]
\begin{center}
\includegraphics[scale=0.5]{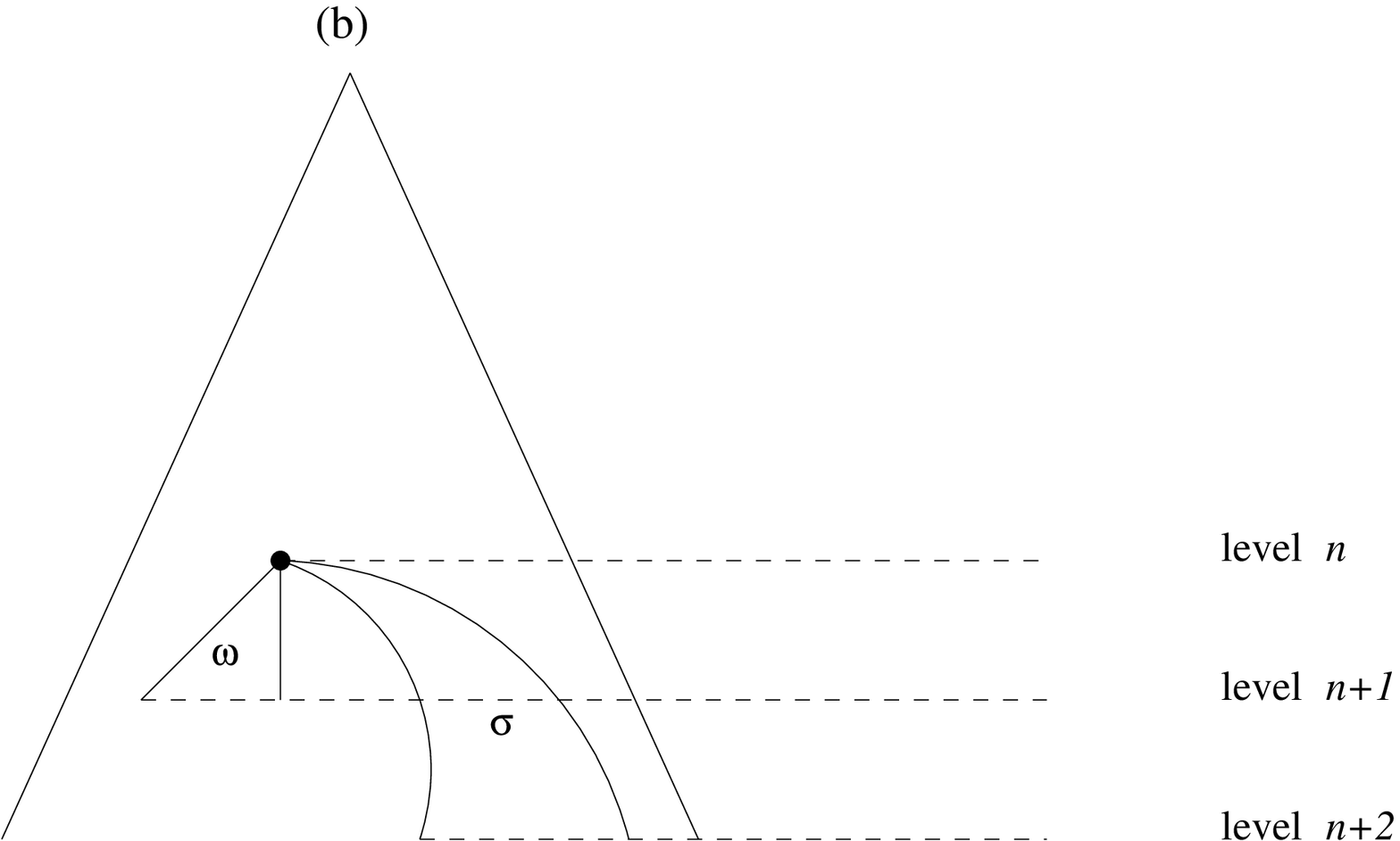}
\end{center}
\end{figure}

As we have declared in the introduction, our main aim is to
determine an expression for the sequence associated with the
doubled rule (with axiom $(b)$) in terms of the sequences
associated with the production rules $\omega$ and $\sigma$. Let's
start by fixing some notations. First of all, $L$ and $M$ will be
the rule operators associated with $\omega$ and $\sigma$,
respectively. Using production rules, our doubled mixed succession
rule will be denoted $(b)\omega^{+1}\sigma^{+2}$, whereas in terms
of rule operators it will be $(b)L^{+1}M^{+2}$. Moreover, we will
denote $p_n (x)$ the level polynomials of
$(b)\omega^{+1}\sigma^{+2}$. Finally, $f^{(b)}(x,t)$ will be the
bivariate generating function of the generating tree, where $t$
keeps track of the level and $x$ keeps track of the label. Our
first result is an expression for $f^{(b)}(x,t)$ in terms of the
rule operators $L$ and $M$.

\begin{prop} Denoting by $^{-1}$ the compositional inverse of an operator, we
have:
\begin{displaymath}
f^{(b)}(x,t)=(\mathbf{1}-\mathbf{t}L-\mathbf{t}^2 M)^{-1}(x^ b).
\end{displaymath}
\end{prop}

\emph{Proof.}\quad The argument to be used here is analogous to
the one used in \cite{FPPR2} for jumping succession rules. Since
each node at level $n$ can be generated either by a node at level
$n-1$ (according to $\omega$) or by a node at level $n-2$
(according to $\sigma$), we have the following expression for $p_n
(x)$:
\begin{displaymath}
p_n (x)=L(p_{n-1}(x))+M(p_{n-2}(x)).
\end{displaymath}

If we impose, by convention, that $p_i (x)=0$, for $i<0$, then the
above expression is meaningful when $n\geq 1$ (recall that, under
our assumptions, $p_0 (x)=x^b$). In order to translate the above
recursion into generating functions, we multiply by $t^n$ both
sides of the above equality and sum up for $n\geq 1$, thus
obtaining:
\begin{displaymath}
\sum_{n\geq 1}p_n (x)t^n =\sum_{n\geq 1}L(p_{n-1}(x))t^n
+\sum_{n\geq 1}M(p_{n-2}(x))t^n ,
\end{displaymath}
whence, using linearity:
\begin{displaymath}
\sum_{n\geq 1}p_n (x)t^n =L\left( \sum_{n\geq 1}p_{n-1}(x)t^n
\right) +M\left( \sum_{n\geq 1}p_{n-2}(x)t^n \right) .
\end{displaymath}

Since $f^{(b)}(x,t)=\sum_{n\geq 0}p_n (x)t^n$, we will then get:
\begin{displaymath}
f^{(b)}(x,t)-x^b =\mathbf{t}L(f^{(b)}(x,t))+\mathbf{t}^2
M(f^{(b)}(x,t)),
\end{displaymath}
whence
\begin{displaymath}
f^{(b)}(x,t)=(\mathbf{1}-\mathbf{t}L-\mathbf{t}^2 M)^{-1}(x^
b).\quad \blacksquare
\end{displaymath}

\bigskip

Expressing the operator $(\mathbf{1}-\mathbf{t}L-\mathbf{t}^2
M)^{-1}$ using power series, we have that
\begin{equation}\label{bivariata}
f^{(b)}(x,t)=\sum_{n\geq 0}\mathbf{t}^n (L+\mathbf{t}M)^n (x^b ).
\end{equation}

Therefore, it is now clear that, if we want to know the sequence
associated with the doubled rule, we need to find an expression
for the binomial $(L+\mathbf{t}M)^n$. In general, this is a
nontrivial problem, since the linear operators $L$ and $M$ usually
\emph{do not commute}. We are then led to first take into
consideration just a special class of pairs of rule operators.
More precisely, in the rest of the paper, we will assume the
following hypothesis:

\bigskip

\begin{center}
\emph{$L$ ed $M$ commutes, i.e. $LM=ML$.}
\end{center}

\bigskip

Using an algebraic terminology, it is said that the
\emph{commutator} $[L,M]=LM-ML$ is equal to zero.

\bigskip

Now let's come back to our problem, that is the determination of
an expression for the binomial $(L+\mathbf{t}M)^n$. We have the
following, crucial result.

\begin{teor}\label{teorema} Denoting by $\mu_r
^{(s)}(x)=\sum_{i}\mu_{r,i}^{(s)}x^i$ the $r$-th level polynomial
of the generating tree of $\sigma$ with axiom $(s)$ and by $(l_n
^{(a)})_{\ninN}$ the numerical sequence associated with $\omega$
with axiom $(a)$, if $(f_n ^{(b)})_{\ninN}$ is the sequence
determined by $(b)\omega^{+1}\sigma^{+2}$, we have:
\begin{equation}\label{main}
f_n ^{(b)} =\sum_{k\geq 0}{n-k\choose
k}\sum_{i}\mu_{k,i}^{(b)}l_{n-2k}^{(i)}.
\end{equation}
\end{teor}

\emph{Proof.} From the fact that $[L,M]=0$ we get immediately:
\begin{displaymath}
(L+\mathbf{t}M)^n =\sum_{k=0}^{n}{n\choose k}\mathbf{t}^k
L^{n-k}M^k .
\end{displaymath}

As a consequence, equality (\ref{bivariata}) can be rewritten as:
\begin{eqnarray*}
f^{(b)}(x,t)&=&\sum_{n\geq 0}\mathbf{t}^n
\sum_{k=0}^{n}\mathbf{t}^k{n\choose k}L^{n-k}M^k (x^b )
\\ &=&\sum_{n\geq 0}\left( \sum_{k\geq 0}{n-k\choose k}L^{n-2k}M^k
(x^b )\right) \mathbf{t}^n .
\end{eqnarray*}

From the above expression we immediately deduce that $p_n
(x)=\sum_{k\geq 0}{n-k\choose k}L^{n-2k}M^k (x^b )$, and so the
$n$-th term of the sequence associated with the doubled rule,
which is $f_n ^{(b)}=p_n (1)$, can be computed as follows:
\begin{eqnarray*}
p_n (1)&=&\left[ \sum_{k\geq 0}{n-k\choose k}L^{n-2k}M^k (x^b
)\right]_{x=1}
\\ &=&\sum_{k\geq 0}{n-k\choose k}\left[ L^{n-2k}(\mu_k
^{(b)}(x))\right]_{x=1}
\\ &=&\sum_{k\geq 0}{n-k\choose
k}\sum_{i}\mu_{k,i}^{(b)}l_{n-2k}^{(i)},
\end{eqnarray*}
and this is precisely our thesis.\cvd

Therefore we have succeeded in finding a formula expressing the
numerical sequence associated with a doubled mixed succession rule
when the related rule operators commute. Specifically, our formula
involves:
\begin{itemize}
\item the distribution of the labels of the production rule
$\sigma$ with axiom $(b)$ inside its generating tree (i.e. the
coefficients $\mu_{k,i}^{(b)}$);

\item the sequences associated with the production rule $\omega$
(i.e. the coefficients $l_{n-2k}^{(i)}$).
\end{itemize}

\section{Examples}

We close by giving two applications of formula (\ref{main}). The
first case is somehow trivial (but leads to interesting
enumerative results), since we deal with the identity operator,
which does not raise any problem concerning commutativity.
However, in any other case, we need to determine some pairs of
commuting rule operators. To this aim, the easiest way is perhaps
to fix a rule operator $L$ and then find the general form of the
rule operators commuting with $L$, which is precisely what we have
done in our second example.

\subsection{The identity operator}

As it is obvious, the identity operator $\mathbf{1}$ commutes with
any linear operator. Therefore, if $L,M$ are any rule operators,
we can consider the two doubled mixed succession rules
$(b)L^{+1}\mathbf{1}^{+2}$ and $(b)\mathbf{1}^{+1}M^{+2}$. Let us
analyze the two cases separately.

\bigskip

Consider first $(b)L^{+1}\mathbf{1}^{+2}$. To apply theorem
\ref{teorema}, we observe that $(l_n ^{(s)})_{\ninN}$ is the
sequence determined by $L$ with axiom $(s)$, whereas $\mu_r
^{(s)}(x)$ is the $r$-th level polynomial of the succession rule
determined by the identity operator $\mathbf{1}$ with axiom $(s)$,
and so it is trivially $\mu_r ^{(s)}(x)=x^s$, whence
\begin{displaymath}
\mu_{r,i}^{(s)}=\left\{ \begin{array}{lll} 1\; ,&\quad&i=s
\\ 0\; ,&\quad&i\neq s
\end{array}\right. .
\end{displaymath}

Thus, denoting by $(f_n ^{(b)})_{\ninN}$ the sequence determined
by the doubled mixed rule $(b)L^{+1}\mathbf{1}^{+2}$, we get:
\begin{displaymath}
f_n ^{(b)}=\sum_{k\geq 0}{n-k\choose k}l_{n-2k}^{(b)}.
\end{displaymath}

Moreover, if $f^{(b)}(x)$ and $l^{(b)}(x)$ are the two generating
functions of the above sequences, standard arguments leads to the
following result.

\begin{cor} For the generating function $f^{(b)}(x)$ we have:
\begin{displaymath}
f^{(b)}(x)=\frac{1}{1-x^2 }\cdot l^{(b)}\left( \frac{x}{1-x^2
}\right) .
\end{displaymath}
\end{cor}

\emph{Examples.}
\begin{enumerate}
\item If $L$ is the rule operator of Catalan numbers described in
(\ref{catalan}), the sequence determined by the mixed rule
$(1)L^{+1}\mathbf{1}^{+2}$ is sequence A105864 in \cite{Sl}, which
has no significant combinatorial interpretation. In order to find
one, we consider a special class of parallelogram polyominoes. A
\emph{1-2 column parallelogram polyomino} is a parallelogram
polyomino whose cells can be either monominoes or dominoes, such
that
\begin{itemize}
\item[i)] each column is entirely made either of monominoes or of
dominoes;

\item[ii)] given a set of consecutive columns starting at the same
height, the leftmost one must be made of monominoes.
\end{itemize}

Such a class of polyominoes can be constructed as follows,
according to the semilength. Given a polyomino $\mathcal{P}$ of
semilength $n$ and such that its rightmost column has $k-1$ cells,
we construct the following set of new polyominoes:
\begin{itemize}
\item either add a new rightmost column made of monominoes ending
at the same height of the rightmost column of $\mathcal{P}$, or

\item add a new cell on the top of the rightmost column of
$\mathcal{P}$ (such a new cell will be a monomino or a domino
according to the type of the column), or

\item add a new rightmost column made of dominoes, starting and
ending at the same heights as the rightmost column of
$\mathcal{P}$.
\end{itemize}

An instance of this ECO construction is shown in figure
\ref{polyomino}.

\begin{figure}[!h]\label{polyomino}
\begin{center}
\includegraphics[scale=0.5]{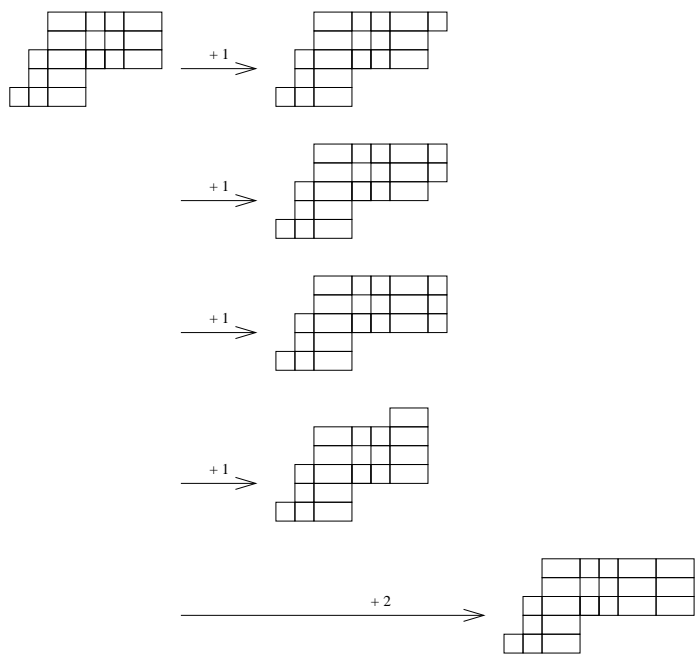}
\end{center}
\caption{Our ECO construction performed on a polyomino of
semilength 14 and such that the rightmost column has 3 cells.}
\end{figure}

As it is clear, in the first two cases polyominoes of semilength
$n+1$ are produced, whereas in the third case a polyomino of
semilength $n+2$ comes out. Now the reader can check that the
above construction can be encoded by the following succession
rule:
\begin{displaymath}
\left\{ \begin{array}{lll} (1)&
\\ (k)&\stackrel{+1}{\rightsquigarrow} (2)(3)\cdots (k)(k+1)
\\ &\stackrel{+2}{\rightsquigarrow} (k)
\end{array}\right. ,
\end{displaymath}
that is precisely the mixed succession rule
$(1)L^{+1}\mathbf{1}^{+2}$.

\item Taking for $L$ the rule operator of Motzkin numbers recalled
in (\ref{motzkin}) (and choosing again $(1)$ as axiom), we get
sequence A128720 of \cite{Sl}. One of the given combinatorial
interpretations for such a sequence is the following: it counts
the number of \emph{2-generalized Motzkin paths}, i.e. paths in
the first quadrant from $(0,0)$ to $(n,0)$ using steps $U=(1,1)$,
$D=(1,-1)$, $h=(1,0)$, and $H=(2,0)$. Various kinds of generalized
Motzkin paths have been extensively studied in the literature, see
for example \cite{dMS,Su}. The mixed succession rule arising in
this case is the following:
\begin{displaymath}
\Omega : \left\{ \begin{array}{lll} (1)&
\\ (k)&\stackrel{+1}{\rightsquigarrow} (1)(2)\cdots (k-1)(k+1)
\\ &\stackrel{+2}{\rightsquigarrow} (k)
\end{array}\right. .
\end{displaymath}

\begin{figure}[!h]\label{motzkin1}
\begin{center}
\includegraphics[scale=0.4]{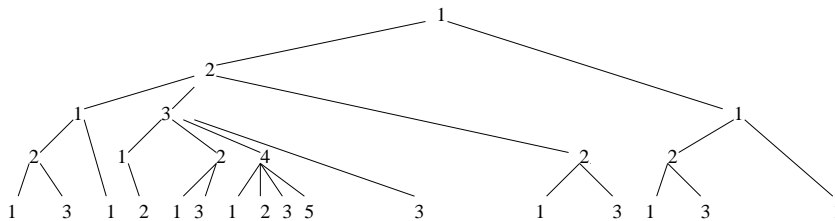}
\end{center}
\caption{The generating tree of $\Omega$}
\end{figure}

\bigskip

It is interesting to notice that $\Omega$ (whose generating tree
is depicted in figure \ref{motzkin1}) indeed describes an ECO
construction for the above class of paths. Leaving the details to
the interested reader, we quickly justify this claim: given a
2-generalized Motzkin path, consider its last descent, i.e. the
final sequence of the path free of $U$ steps. Construct a set of
new paths as follows: either replace each $h$ step with a $U$ step
and add a $D$ step at the end, or just add a $h$ step at the end,
or simply add a $H$ step at the end. In the first two cases the
length of the path is increased by 1, whereas in the last case it
is increased by 2. It is now easy to show that this construction
is encoded precisely by the mixed succession rule $\Omega$. We
also remark that another interpretation of sequence A128720 is
provided in \cite{Sl}, namely using skew Dyck path. It would be
interesting to use the above mixed rule to describe a construction
for this latter combinatorial structure as well.

To conclude this example, we also notice that, according to
\cite{BMS}, the rule $\Omega$ is alternatively described by the
$A$-matrix
\begin{displaymath}
\left( \begin{array}{cccccc} 0&1&0&0&0&\cdots
\\ 1&0&1&1&1&\cdots
\end{array}\right) ,
\end{displaymath}
which implies that the associated ECO matrix is actually a Riordan
array. Its $A$-sequence has generating function:
\begin{eqnarray*}
A(t)&=&\frac{1-t+t^2 +\sqrt{1-2t+7t^2 -10t^3 +5t^4 }}{2(1-t)}
\\&=&1+
2t^2 +t^3 -t^4 +6t^6 +5t^7 -16t^8 +\cdots ,
\end{eqnarray*}
and this shows that a direct dependence of row $n+1$ from row $n$
is very unlikely. Using the theory of Riordan arrays, we can
determine the formal power series $d(t)$ and $h(t)$ defining our
ECO matrix. More precisely, since $h(t)=tA(h(t))$, we find:
\begin{eqnarray*}
h(t)&=&\frac{1}{2}\left( 1-\sqrt{\frac{1-3t-t^2 }{1+t-t^2
}}\right)
\\&=&t+2t^3 +t^4 +7t^5 +10t^6 +37t^7 +82t^8 +\cdots .
\end{eqnarray*}

Since column 0 is not privileged, we have $d(t)=\frac{h(t)}{t}$,
so our Riordan array is completely determined. Denoting by
$S(t)=\sum_{n\geq 0}S_n t^n$ the generating function of the row
sums of the array, since $S(t)=\frac{d(t)}{1-h(t)}$, we get:
\begin{eqnarray*}
S(t)&=&\frac{1-t-t^2 -\sqrt{(1-t-t^2 )(1-3t-t^2 )}}{2t^2 }
\\&=&1+t+3t^2 +6t^3 +16t^4 +40t^5 +109t^6 +297t^8 +\cdots.
\end{eqnarray*}

Using standard methods of asymptotic analysis, we obtain the
asymptotic value:
\begin{displaymath}
S_n \sim \frac{K}{2}\cdot
\frac{\psi^{-n-2}}{4^{n+2}(2n+3)}{2n+4\choose n+2},
\end{displaymath}
where $K=\sqrt{22\sqrt{13}-78}\approx 1.149838276...$ and $\psi
=\frac{\sqrt{13}-3}{2}$. For $n=60$ it is $S_{60}\approx
4.960467337...\times 10^{28}$, whereas the above formula gives
$4.949459297...\times 10^{28}$, with a relative error of 0.22\%.

\item If $L$ is as in (\ref{bell}), defining Bell numbers, the
resulting sequence $(f_n ^{(1)})_{\ninN}$ starts
$1,1,3,7,22,75,\ldots$ and is not recorded in \cite{Sl}. Thanks to
our theory, it is possible to give a combinatorial interpretation
to such a sequence, by performing an ECO construction described by
the mixed succession rule:
\begin{displaymath}
\Omega : \left\{ \begin{array}{lll} (1)&
\\ (k)&\stackrel{+1}{\rightsquigarrow} (k)^{k-1} (k+1)
\\ &\stackrel{+2}{\rightsquigarrow} (k)
\end{array}\right. .
\end{displaymath}

We call a \emph{lacunary partition} of $[n]=\{ 1,2,\ldots ,n\}$
any partition of a subset $S$ of $[n]$ such that $[n]\setminus S$
is disjoint union of intervals of even cardinality. For instance,
the partition $\{ \{ 1,8,12\} ,\{ 2\} ,\{ 3,9\} \}$ is a lacunary
partition of $[14]$. An ECO construction for the class of lacunary
partitions works as follows: given a lacunary partition $\pi$ of
$[n]$, construct a set of new lacunary partitions by either adding
the block $\{ n+1\}$, or adding $n+1$ to each of the block of
$\pi$, or else leaving $\pi$ unchanged, but thinking of it as a
lacunary partition of $[n+2]$. Of course, performing one of the
first two operations leads to a lacunary partition of $[n+1]$,
whereas the last one produces a lacunary partition of $[n+2]$. The
reader can now check that such a construction is encoded by
$\Omega$.
\end{enumerate}

Now consider $(b)\mathbf{1}^{+1}M^{+2}$. In this case, the
sequence $(l_n ^{(s)})_{\ninN}$ is the one determined by the
identity operator, and so $l_n ^{(s)}=1$, for all $\ninN$. On the
other hand, the polynomial $\mu_r ^{(s)}(x)$ is the level
polynomial of the rule associated with $M$ with axiom $(s)$.
Applying theorem \ref{teorema}, for the sequence $(f_n
^{(b)})_{\ninN}$ determined by $(b)\mathbf{1}^{+1}M^{+2}$ we then
get:
\begin{displaymath}
f_n ^{(b)}=\sum_{k\geq 0}{n-k\choose
k}\sum_{i}\mu_{k,i}^{(b)}=\sum_{k\geq 0}{n-k\choose k}m_k ^{(b)},
\end{displaymath}
where, of course, $(m_n ^{(s)})_{\ninN}$ is the sequence
associated with the rule operator $M$ with axiom $(s)$.

As in the preceding case, we can easily recover the generating
function of the sequence $f_n ^{(b)}$ starting from that of the
sequence $m_n ^{(b)}$.

\begin{cor} For the generating function $f^{(b)}(x)$ we have:
\begin{displaymath}
f^{(b)}(x)=\frac{1}{1-x}\cdot m^{(b)}\left( \frac{x^2
}{1-x}\right) .
\end{displaymath}
\end{cor}

\emph{Examples.} We leave as an open problem that of finding an
ECO construction described by $(1)\mathbf{1}^{+1}M^{+2}$ for each
of the structures mentioned in the next examples. In the first of
them, we also provide a rigorous proof (via Riordan arrays) of the
fact that a certain sequence on \cite{Sl} comes out, whereas the
details of the second example are left to the reader.
\begin{enumerate}
\item If $M$ is the rule operator of Catalan numbers described in
(\ref{catalan}), the sequence determined by the mixed rule
$(1)\mathbf{1}^{+1}M^{+2}$ is sequence A090344 in \cite{Sl}. Such
a sequence counts Motzkin paths without horizontal steps at odd
height.

From an enumerative point of view, observe that
$(1)\mathbf{1}^{+1}M^{+2}$ can be explicitly written as:
\begin{displaymath}
\left\{ \begin{array}{lll} (1)&
\\ (k)&\stackrel{+1}{\rightsquigarrow} (k)
\\ &\stackrel{+2}{\rightsquigarrow} (2)(3)\cdots (k)(k+1)
\end{array}\right. .
\end{displaymath}

According to \cite{BMS}, this corresponds to the $A$-matrix
\begin{displaymath}
\left( \begin{array}{cccccc} 1&1&1&1&1&\cdots
\\ 0&1&0&0&0&\cdots
\end{array}\right) ,
\end{displaymath}
which gives rise to the following \emph{vertically stretched
Riordan array} (see \cite{CMS}):
\begin{displaymath}
\begin{array}{c|ccccc}
n\backslash k & 1 & 2 & 3 & 4 & 5 \\ \hline
0 & 1 &  &  &  &  \\
1 & 1 &  &  &  &  \\
2 & 1 & 1 &  &  &  \\
3 & 1 & 2 &  &  &  \\
4 & 1 & 4 & 1 &  &  \\
5 & 1 & 7 & 3 &  &  \\
6 & 1 & 13 & 8 & 1 &  \\
7 & 1 & 24 & 18 & 4 &  \\
8 & 1 & 47 & 40 & 13 & 1 \\
\end{array}
\end{displaymath}

The entries $m_{n,k}$ of the above array obey the recursion
$m_{n+1,k+1}=m_{n,k+1}+\sum_{j=0}^{\infty}m_{n-1,k+j}$. To obtain
a proper Riordan array (to be denoted $(d(t),h(t))$), we can
simply shift column $k+1$ up by $k$ position (for all $k$), so
that the entries $p_{n,k}$ of the new (lower triangular) array
will now satisfy the equalities
$p_{n+1,k+1}=p_{n,k+1}+\sum_{j=0}^{\infty}p_{n-j,k+j}$. These
recurrence relations are described by the (infinite) $A$-matrix:
\begin{displaymath}
\left( \begin{array}{cccccc}
\vdots & \vdots & \vdots & \vdots & \vdots &{\rotatebox{90}{$\ddots$}} \\
0 & 0 & 0 & 0 & 1 & \cdots \\
0 & 0 & 0 & 1 & 0 & \cdots \\
0 & 0 & 1 & 0 & 0 & \cdots \\
0 & 1 & 0 & 0 & 0 & \cdots \\
1 & 1 & 0 & 0 & 0 & \cdots \\
\end{array} \right).
\end{displaymath}

To get the generating function of the $A$-sequence of our Riordan
array we can use the equality $A(t)=\sum_{j=0}^{\infty}t^j
A(t)^{-j}P^{[j]}(t)$ (shown in \cite{BMS}), where $P^{[j]}(t)$ is
the generating function of row $j$ in the $A$-matrix. In our case,
we have $P^{[0]}(t)=1+t$ and $P^{[j]}(t)=t^j$, for $j>0$, which
yields:
\begin{eqnarray*}
A(t)&=&1+t+tA(t)^{-1}t+t^2 A(t)^{-2}t^2 +t^3 A(t)^{-3}t^3 +\cdots
\\ &=&t+1+\frac{t^2}{A(t)}+\frac{t^4}{A(t)^2}+\frac{t^6}{A(t)^3}+\cdots
=t+\frac{A(t)}{A-t^2}.
\end{eqnarray*}

The solution of the above equation is:
\begin{eqnarray*}
A(t)&=&\frac{1+t+t^2 +\sqrt{1+2t+3t^2 -2t^3 +t^4}}{2}
\\&=&1+t+t^2 -t^3 +t^4 -2t^6 +4t^7 -3t^8 -5t^9 +\cdots.
\end{eqnarray*}

Now we are ready to compute $d(t)$ and $h(t)$. Using the well
known formula $h(t)=tA(h(t))$, we get an equation of degree two,
whose solution is:
\begin{eqnarray*}
h(t)&=&\frac{1}{2t}\left( 1-\sqrt{\frac{1-t-4t^2}{1-t}}\right)
\\ &=&t+t^2 +2t^3 +3t^4 +6t^5 +11t^6 +23t^7 +47t^8 +102t^9 +\cdots
,
\end{eqnarray*}
and obviously $d(t)=\frac{1}{1-t}$. Our original (stretched)
Riordan array is therefore:
\begin{displaymath}
\left( \frac{1}{1-t},\ \frac{1}{2}\left(
1-\sqrt{\frac{1-t-4t^2}{1-t}}\right) \right) ,
\end{displaymath}
and the generating function of the sequence of its row sums is:
\begin{displaymath}
S(t)=\sum_{n=0}^{\infty}S_n t^n
=\frac{d(t)}{1-h(t)}=\frac{1}{2t^2}\left(
1-\sqrt{\frac{1-t-4t^2}{1-t}}\right) .
\end{displaymath}

It is possible to find an asymptotic value for $S_n$. The equation
$1-t-4t^2 =0$ has the two solutions $r_1
=\frac{-\sqrt{17}+1}{2}\approx -0.6403882032...$ and $r_2
=\frac{\sqrt{17}-1}{2}\approx 0.3903882032...$, and therefore this
last value is the dominating singularity. Using Bender's theorem
\cite{B}, we obtain:
\begin{eqnarray*}
S_n &=&[t^n ]\frac{1}{2t^2}\left( 1-\sqrt{\frac{1-t-4t^2}{1-t}}
\right)
\\ &\approx&-\frac{1}{2}\left[ \left[ \sqrt{\frac{1-t/r_1}{1-t}}\ \right]_{t=r_2}\right] [t^{n+2}]\sqrt{1-\frac{t}{r_2}}
\\ &=&\frac{K}{2}\frac{1}{(2n+3)r_2^{n+2}4^{n+2}}{2n+4\choose n+2}.
\end{eqnarray*}

In the above formula, we have used the notation $[t^n ]f(t)$ to
denote the coefficient of $t^n$ in the formal power series $f(t)$.
For instance, for $n=50$ we have $S_{50}\approx
1.091877333...\times 10^{18}$, while the approximate value is
$1.075272279...\times 10^{18}$, with a relative error of 1.54\%.
Since $r_1 <1$, its contribution decreases as $n$ increases, but
for small values of $n$ it cannot be ignored.

\item Taking for $M$ the rule operator of Motzkin numbers recalled
in (\ref{motzkin}) (and choosing again $(1)$ as axiom), we get
sequence A026418 of \cite{Sl}. It counts ordered trees having no
branches of length 1, according to the number of edges.

\item If $M$ is as in (\ref{bell}), defining Bell numbers, the
resulting sequence $(f_n ^{(1)})_{\ninN}$ starts
$1,1,2,3,6,11,23,47,103,\ldots$ and is not recorded in \cite{Sl}.
\end{enumerate}

\subsection{A factorial-like rule operator}

Consider the rule operator $L=\mathbf{x}^2 D$ associated with the
production rule
\begin{displaymath}
\omega :\quad (k)\rightsquigarrow (k)^{k+1}.
\end{displaymath}

We start by determining the family of sequences related to $L$.

\begin{lemma}\label{sequenza} If $(l_n ^{(b)})_{\ninN}$ is the sequence determined
by $L$ with axiom $(b)$, then we have, for all $\ninN$:
\begin{displaymath}
l_n ^{(b)}=(n+b-1)_{b-1}=\frac{(n+b-1)!}{(b-1)!}=n!{n+b-1\choose
b-1},
\end{displaymath}
where $(x)_y=x(x-1)\cdot \ldots \cdot (x-y+1)$ denotes the usual
falling factorial.
\end{lemma}

\emph{Proof (sketch).}\quad Use a simple induction argument. For
$b=1$ it is well-known (see, for instance, \cite{FP1}) that $l_n
^{(1)}=n!$. Now observe that the recursion defined by the
associated production rule implies that $l_{n+1}^{(b)}=bl_n
^{(b+1)}$, whence it is easy to derive the thesis.\cvd

According to our program, we start by computing the general form
of a rule operator commuting with $L$.
\begin{teor}
Let $M$ be a rule operator such that $M(1)=x^a$, for some $a\in
\mathbf{N}$. Then, $M$ commutes with $L$ if and only if
\begin{displaymath}
M=L_{[a]}=\frac{\mathbf{x}^{a+1}}{(a-1)!}D^a \mathbf{x}^{a-1}.
\end{displaymath}
\end{teor}

\emph{Proof.}\quad Suppose that $M$ commutes with $L$. On the
polynomial $1$ it is
\begin{displaymath}
M(x)=M(\mathbf{x}^2 D(1))=\mathbf{x}^2 D(M(1))=\mathbf{x}^2 D(x^a)
= a x^{a+1}.
\end{displaymath}

Now suppose by induction that $M(x^n )=n{a+n-1\choose
a-1}x^{a+n}$. Then, on $x^{n+1}$ we have:
\begin{eqnarray*}
M(x^{n+1})&=&\frac{1}{n}M(nx^{n+1})=\frac{1}{n}M(\mathbf{x}^2
D(x^n ))=\frac{1}{n}\mathbf{x}^2 D(M(x^n ))
\\ &=&\frac{1}{n}\mathbf{x}^2 D\left( n{a+n-1\choose
a-1}x^{a+n}\right) =(a+n){a+n-1\choose a-1}x^{a+n+1}
\\ &=&(n+1){a+n\choose a-1}x^{a+n+1}.
\end{eqnarray*}

We have thus showed that $M(x^n )=n{a+n-1\choose a-1}x^{a+n}$,
that is $M=L_{[a]}
=\frac{\mathbf{x}^{a+1}}{(a-1)!}D^a\mathbf{x}^{a-1}$, as
desired.

As far as the converse is concerned, we leave to the reader the
proof of the fact that the operator
$L_{[a]}=\frac{\mathbf{x}^{a+1}}{(a-1)!}D^a \mathbf{x}^{a-1}$
commutes with $L$.\cvd

Consider now the case $a=2$, so to obtain the operator $M=L_{[2]}
=\mathbf{x}^3 D\mathbf{x}^2$. The mixed succession rule
$\Omega_{b}=(b)L^{+1}M^{+2}$ is
\begin{equation}\label{permut}
\Omega_b : \left\{ \begin{array}{ll} (b)&
\\ (k)&\stackrel{+1}{\rightsquigarrow} (k+1)^k
\\ &\stackrel{+2}{\rightsquigarrow} (k+2)^{k(k+1)}
\end{array}\right. .
\end{equation}

The first levels of the corresponding generating tree, when $b=1$,
appear as in figure 3.

\begin{figure}[!h]\label{alberofattoriale}
\begin{center}
\includegraphics[scale=0.5]{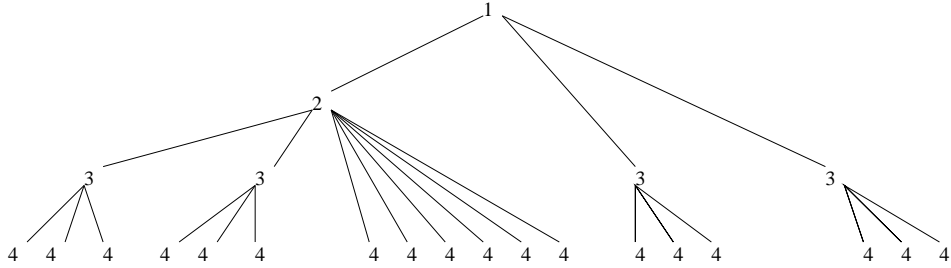}
\end{center}
\caption{The first levels of the generating tree of
$(1)L^{+1}M^{+2}$.}
\end{figure}

In order to apply theorem \ref{teorema}, we need to know the
sequence $l_n ^{(s)}$ determined by the production rule associated
with $L$ with axiom $(s)$ and the level polynomials $\mu_r
^{(s)}(x)$ of the generating tree related to $M$ with axiom $(s)$.

The first information is provided by lemma \ref{sequenza}, that is
$l_n ^{(s)}=n!{n+s-1\choose s-1}$. As far as the polynomials
$\mu_r ^{(s)}(x)$ are concerned, we observe that, in the
generating tree associated with $M$, only one label appears at any
given level; more precisely, the only label at level $r$ is
$(s+2r)$. Therefore $\mu_r ^{(s)}(x)$ consists of only one
monomial, and we simply have to determine its coefficient. The
following, simple lemma finds this coefficient.

\begin{lemma} The generating tree associated with $M$ having axiom
$(s)$ has $(s)^{2r}$ nodes at level $r$ (each of which is labelled
$(s+2r)$), where $(x)^{y}=x(x+1)\cdot \ldots \cdot (x+y-1)$
denotes the raising factorial.
\end{lemma}

\emph{Proof.}\quad At level $0$ and $1$ there are, respectively, 1
and $s(s+1)$ nodes. By induction, suppose that at level $r$ we
have $(s)^{2r}$ nodes; since each of them is labelled $(s+2r)$, it
produces $(s+2r)(s+2r+1)$ sons, whence the thesis immediately
follows.\cvd

As a consequence, we have that $\mu_r ^{(s)}(x)=(s)^{2r}x^{s+2r}$,
which means that $\mu_{r,s+2r}^{(s)}=(s)^{2r}$, whereas
$\mu_{r,j}^{(s)}=0$, for $j\neq s+2i$.

We are now ready to apply theorem \ref{teorema}, thus getting for
the sequence $(f_n ^{(b)})_{\ninN}$ determined by
$(b)L^{+1}M^{+2}$ the following formula:
\begin{eqnarray}\label{fact}
f_n ^{(b)}&=&\sum_{k\geq 0}{n-k\choose
k}(b)^{2k}(n-2k)!{n+b-1\choose b+2k-1}\nonumber
\\ &=&\frac{(n+b-1)!}{(b-1)!}\sum_{k\geq 0}{n-k\choose k}\nonumber
\\ &=&(b)^n F_n ,
\end{eqnarray}
where $(F_n )_{\ninN}$ is the sequence of Fibonacci numbers.

To give a combinatorial interpretation for the sequence $(f_n
^{(b)})_{\ninN}$ we refer to the generating tree of the mixed rule
$\Omega_b$. In what follows, $\overline{S}_n$ will denote the set
of \emph{coloured permutations} of $[n]$, i.e. permutations whose
elements can possibly be coloured (a coloured element will simply
be overlined). Moreover, we introduce the notion of \emph{paired
coloured permutation}, to mean a coloured permutation such that,
if $I$ denotes the set of coloured elements, then $I$ is a
disjoint union of intervals of even cardinality. So, for instance,
the permutation
$\overline{93}1\overline{4}6\overline{5}7\overline{28}$ is a
paired coloured permutation belonging to $\overline{S}_9$.

\begin{prop} Given $b\in \mathbf{N}$, fix a permutation $\tau \in
S_{b-1}$. Then $f_n ^{(b)}$ is the number of paired coloured
permutations $\pi \in \overline{S}_{b+n-1}$ in which the elements
$1,2,\ldots ,b-1$ are not coloured and appear in $\pi$ as a
pattern isomorphic to $\tau$.
\end{prop}

\emph{Proof.}\quad Let $\pi \in \overline{S}_{b+k-1}$ be a paired
coloured permutation satisfying the hypotheses of the theorem.
Starting from $\pi$ we construct a new sets of permutations as
follows:
\begin{enumerate}
\item add the noncoloured element $b+k$ in any of the $b+k$
possible positions, so to obtain $b+k$ new permutations belonging
to $\overline{S}_{b+k}$; \item add the two coloured elements
$\overline{b+k}$ and $\overline{b+k+1}$ in any possible positions:
this can be done in $(b+k)(b+k+1)$ different ways, and produces
permutations belonging to $\overline{S}_{b+k+1}$.
\end{enumerate}

Moreover, the additional hypothesis that the subpermutation of
$\pi$ constituted by the elements $1,2,\ldots ,b-1$ must be
isomorphic to $\tau$ implies that $\tau$ is the minimal
permutation of our class, and any other permutation obtained using
the above described construction must avoid the pattern $\tau$ in
its elements $1,2,\ldots ,b-1$. It is now easy to recognize that,
if a permutation $\pi \in \overline{S}_{b+k-1}$ is given label
$(b+k)$, the above construction can be described by the mixed
succession rule $\Omega$ given in (\ref{permut}), which is enough
to conclude.\cvd

\begin{cor} If $b=1$, then $f_n =f_n ^{(1)}=n!F_n$, and the sequence $(f_n
)_{\ninN}$ enumerates the class of paired coloured permutations.
\end{cor}

\emph{Proof.}\quad Set $b=1$ in the previous proposition.\cvd

This last sequence also appears in \cite{Sl} (it is essentially
sequence A005442), and can be obtained as the row sums of a
particular convolution matrix (see \cite{K}). The combinatorial
interpretation of $(f_n )_{\ninN}$ reported in \cite{Sl} seems to
be essentially different from the one given here: it would be
interesting to have a bijective argument explaining how to relate
these two interpretations.

\section{Final remarks}

In the present paper we have studied doubled mixed succession
rules, and, in the commutative case, we have been able to give an
expression for the sequence associated with one of such rules in
terms of the sequences associated with the constituent simple
succession rules. The next step should be to have an analogous
result for more general kinds of doubled mixed succession rules.
For instance, one could consider two succession rules whose
associated rule operators obey some weaker form of commutativity,
such as $LM=qML$, for a given scalar $q$ (or, more generally,
$LM=f(q)ML$, for some polynomial $f$).

Another presumably fertile line of research concerns exhaustive
generation. Similarly to what has been done for classical
succession rules, one can try to develop general exhaustive
generation algorithms based on mixed succession rules, maybe
finding a new way of defining general Gray codes depending only on
the form of the mixed succession rule under consideration.

\end{document}